\documentclass[10pt,twoside]{article}
\usepackage{graphicx}
\usepackage{amsmath}
\usepackage{Latex-document}

\title{\bf  Free Probability and Combinatorics\vskip 6mm}
\author{P. Biane\vspace*{-0.5cm}\thanks{D\'epartement de Math\'ematiques et
Applications, \'Ecole Normale Sup\'erieure, 45 rue d'Ulm 75005
Paris, France. E-mail: Philippe.Biane@ens.fr}}
\date{\vspace{-8mm}}

\begin{document}

\maketitle

\thispagestyle{first} \setcounter{page}{765}

\begin{abstract}

\vskip 3mm

A combinatorial approach to free probability theory has been developped by
Roland Speicher, based on the notion of noncrossing cumulants,
a free  analogue of the classical theory of cumulants in
probability theory. We review this theory, and explain
 the connections between free probability theory
and random matrices. We relate  noncrossing cumulants to classical
cumulants and  also to characters of large symmetric groups.
Finally we give applications to the asymptotics of representations
of symmetric groups, specifically to the Littlewood-Richardson
rule.

\vskip 4.5mm


\end{abstract}

\vskip 12mm

\section{Introduction} \label{section 1}\setzero

\vskip-5mm \hspace{5mm}

Free probability has been introduced by D. Voiculescu \cite{V1}
as a means of studying the
group von Neumann algebras of free groups, using probabilistic techniques.
His theory has become very successful when he discovered
a deep relation with the theory of random
matrices, and solved some old questions in operator algebra, see
\cite{B4}, \cite{G}, \cite{V4} for an overview.
 A purely combinatorial approach to Voiculescu's definition of freeness
has been given by R. Speicher \cite{S1}, \cite{S2},
building on G. C. Rota's \cite{R} approach to classical
probability. It is based on the notion of noncrossing partitions, also known as
``planar diagrams'' in quantum field theory, and provides unifying concepts for
many computations in free probability.
Noncrossing partitions turn out to be connected with the geometry of the
symmetric group, and this leads to some new understanding of the asymptotic
behaviour of the characters and
representations of large symmetric groups. Our aim is to
survey these results, we shall start with the basic definition of freeness, then
explain its connection to random matrix theory. In the third section
we review Speicher's theory. In the fourth section we show how noncrossing
 cumulants
arise naturally in connection with classical cumulants associated
with random matrices,  and with
characters of symmetric groups. Finally in section 5 we explain the asymptotic
behaviour of representations of symmetric groups in terms of free probability
concepts.

\section{Freeness and random matrices} \label{section 2}

\setzero\vskip-5mm \hspace{5mm}

The usual framework for free probability is a von Neumann algebra $A$, equipped
with a faithful, tracial,  normal state $\tau$.
To any self-adjoint element $X\in A$ one can associate its distribution,
 the probability measure on
the real line, uniquely determined by the identity
$\tau(X^n)=\int_{\bf R}x^n\mu(dx)$ for all $n\geq 1$. This makes it natural to
think of the elements of $A$ as noncommutative random variables, and of
$\tau$ as an expectation map, and
 one usually calls noncommutative probability space
such a pair $(A,\tau)$.
Although a great deal of the
theory, especially the combinatorial side, can be developped in a purely
algebraic way, assuming only that $A$ is a complex algebra with unit, and $\tau$
a complex linear functional,
we shall stick to the von Neumann framework in the present exposition.

Given $(A,\tau)$, one considers a family $\{A_i;i\in I\}$ of von Neumann
subalgebras. This family is called a {\it free family}
 if the following holds: for any
$k\geq 1$ and $k$-tuple
$a_1,\ldots, a_k\in A$ such that
\begin{itemize}
\item each $a_j$ belongs to some
algebra $A_{i_j}$,
with
$i_1\ne i_2,i_2\ne i_3,\ldots ,i_{k-1}\ne i_{k}$,
\item $\tau(a_j)=0$
 for all $j$,
 \end{itemize}  one has $\tau(a_1\ldots
a_k)=0$.

Moreover, a
family of elements of $A$ is called free if the von Neumann algebras
each of them generates form a free family.
  Freeness is a
noncommutative notion analogous to the  independence of $\sigma$-fields
 in
probability theory, but which incorporates also the  notion of algebraic
independence.

Observe that if $a_1$ and $a_2$ are free elements in $(A,\tau)$, and one defines
the centered elements $\hat a_i=a_i-\tau(a_i)1$ then one can
conpute
$$\tau(a_1a_2)=\tau(\hat a_1\hat a_2)+\tau(a_1)\tau(a_2)=\tau(a_1)\tau(a_2)$$
where the freeness condition has been used to get $\tau(\hat a_1\hat a_2)=0$.
 Actually, if $\{A_i;i\in I\}$ is a free family,
  it is not dificult to see that one can compute the
 value of $\tau$ on any product of the form $a_1\ldots a_k$,
 where each $a_j$ belongs to
 some of the $A_i\, 's$, in terms of the quantities
 $\tau(a_{j_1}\ldots a_{j_l})$ where all the elements $a_{j_1},\ldots, a_{j_l}$
  belong to the same
 subalgebra.  This implies that the value of $\tau$ on the algebra generated by
 the family $\{A_i;i\in I\}$
  is completely determined by the restrictions of $\tau$ to each
 of these subalgebras. However the problem of finding an explicit formula is
 nontrivial, and this is where combinatorics comes in.
 We shall describe  Speicher's theory
 of noncrossing cumulants, which solves this
 problem, in the next section, but
 before that we explain how free probability is relevant to understand
 large  random matrices.

 Consider $n$ random $N\times N$ matrices $X_1^{(N)},\ldots ,X_n^{(N)}$,
 of the form
 \begin{equation}\label{rm}
 X_j^{(N)}=U_jD_j^{(N)}U_j^*
 \end{equation}
 where $D_j^{(N)}; j=1,\ldots,n$
 are diagonal, hermitian, nonrandom
 matrices and $U_j$
 are independent
 unitary random matrices, each distributed with the Haar measure on
 the unitary group ${\bf U}(N)$.
 In other words we have fixed the spectra of the $X_i^{(N)}$ but
 their eigenvectors
 are chosen at random. The  $n$-tuple $X_1^{(N)},\ldots, X_n^{(N)}$ can be
 recovered, up to a global unitary
 conjugation $X_i^{(N)}\mapsto UX_i^{(N)}U^*$, (where $U$ does not depend on
 $i$), from its mixed moments, i.e. the set of complex numbers
 $\frac{1}{N}Tr(X_{i_1}^{(N)}\ldots X_{i_k}^{(N)})$ where
  $i_1,\ldots i_k$ are arbitrary sequences of indices in $\{1,\ldots,n\}$.
  In particular the spectrum of any
 noncommutative polynomial of the $X_i^{(N)}$ can be recovered from these data.
 A most remarkable fact is that if we assume that the individual moments
 $\frac{1}{N}Tr((X_i^{(N)})^k)$ converge as $N$ tends to infinity,
 then the mixed moments
 $\frac{1}{N}Tr(X_{i_1}^{(N)}\ldots X_{i_k}^{(N)})$ converge
  in probability, and their limit
 is obtained by the prescriptions of free probability.

 {\bf Theorem 1.} \it Let
 $(A,\tau)$ be a noncommutative probability space
 with free selfadjoint elements
 $X_1,\ldots, X_n$, satisfying $\tau(X_i^k)
 =\lim_{N\to\infty}\frac{1}{N}Tr((X_i^{(N)})^k)$, for all $i$ and $k$,
  then, in probability,
 $\frac{1}{N}Tr(X_{i_1}^{(N)}\ldots X_{i_k}^{(N)})\to_{N\to\infty}
 \tau(X_{i_1}\ldots X_{i_k})$, for all $i_1,\ldots ,i_k$.  \rm

 This striking result was first proved by D. Voiculescu \cite{V3}, and has lead to the
 resolution of many open problems about von Neumann algebras, upon which we shall
 not touch here.

\section{Noncrossing partitions and  cumulants} \label{section 3}

\setzero\vskip-5mm \hspace{5mm}

A partition of the set $\{1,\ldots, n\}$ is said to have a crossing
if there exists
a quadruple $(i,j,k,l)$, with $1\leq i<j<k<l\leq n$, such that $i$
and
$k$  belong to some class of the partition and  $j$ and $l$
 belong to another class. If a partition has no crossing, it is called
noncrossing. The set of all noncrossing partitions of
$\{1,\ldots, n\}$ is denoted by
$NC(n)$. It is a lattice for the refinement order, which
 seems to have been first systematically investigated in \cite{K}.

 Let
$( A,\tau)$
 be a non-commutative probability space, then
 we shall define a family  $R^{(n)}$ of
 $n$-multilinear
 forms on $ A$, for $n\geq 1$,  by the following formula
 \begin{equation}\label{NCN}
 \tau(a_1\ldots a_n)=\sum_{\pi\in NC(n)}R[\pi](a_1,\ldots ,a_n).
 \end{equation}
 Here, for $\pi\in NC(n)$, one has defined
 $$R[\pi](a_1,\ldots ,a_n)=\prod_{V\in
\pi}R^{(\vert V\vert)}(a_V)$$
where  $a_V=(a_{j_1},\ldots, a_{j_k})$  if $V=\{j_1,\ldots ,j_k\}$ is a
class of the partition $\pi$, with
$j_1<j_2<\ldots <j_k$ and $|V|=k$ is the number of elements of $V$.
  In particular $R[1_n]=R^{(n)}$ if $1_n$ is
the
partition with only one class.
Thus one has, for $n=3$,
$$
\begin{array}{lcl}
\tau(a_1a_2a_3)&=&R^{(3)}(a_1,a_2,a_3)+R^{(2)}(a_1,a_2)R^{(1)}(a_3)
+R^{(2)}(a_1,a_3)R^{(1)}(a_2)\\ &&+R^{(2)}(a_2,a_3)R^{(1)}(a_1)+
R^{(1)}(a_1)R^{(1)}(a_2)R^{(1)}(a_3).
\end{array}
$$

 Observe that
 $$\tau(a_1\ldots a_n)=R^{(n)}(a_1,\ldots,a_n)+\text{terms involving
 $R^{(k)}$ for $k<n$}$$ so that the $R^{(n)}$ are well
  defined by (\ref{NCN}) and can be computed by induction on
 $n$. They are called the noncrossing (or sometimes free)
  cumulant functionals on $A$.

 The formula (\ref{NCN}) can be inverted to yield
 $$R^{(n)}(a_1,\ldots, a_n)=\sum_{\pi\in NC(n)}Moeb([\pi,1_n])\tau
 [\pi](a_1,\ldots ,a_n).$$
 Here $\tau[\pi](a_1,\ldots,a_n)=
 \prod_{V\in \pi}\tau(a_{j_1}\ldots a_{j_k})$ where
  $V=\{j_1,\ldots ,j_k\}$  are the classes of $\pi$, and $Moeb$ is the M\"obius
  function of the lattice $NC(n)$, see \cite{S2}.

 For example, one has
 $$\begin{array}{cl}
 R^{(1)}(a_1)=\tau(a_1);&\qquad
 R^{(2)}(a_1,a_2)=\tau(a_1a_2)-\tau(a_1)\tau(a_2);\\
R^{(3)}(a_1,a_2,a_3)=&\tau(a_1a_2a_3)
-\tau(a_1)\tau(a_2a_3)-\tau(a_2)\tau(a_1a_3)
\\&-\tau(a_3)\tau(a_1a_2)+2\tau(a_1)\tau(a_2)\tau(a_3).
\end{array}
$$
 Note that when the lattice of all partitions is used instead of noncrossing
 partitions, then one gets the usual family of cumulants (see Rota \cite{R}),
  with
 another M\"obius function.

  The connection
between noncrossing cumulants and freeness is the following result from
section 4 of \cite{S1}.

{\bf Theorem 2.} \it Let $\{A_i;i\in
I\}$ be a free family of subalgebras of
$( A,\tau)$, and $a_1,\ldots,
 a_n\in  A$ be such that  $a_j$ belongs to some $A_{i_j}$
for each
$j\in\{1,2,\ldots ,n\}$. Then one has
$R^{(n)}(a_1,\ldots, a_n)=0$ if there exists some
$j$
and
$k$ with
$i_j\not=i_k$.\rm

This result leads to an explicit expression for
$\tau(a_1\ldots a_n)$, where
 $a_1,\ldots,a_n$ is an arbitrary sequence in $
A$,
such that each $a_j$ belongs to one of the algebras $A_i;i\in I$.
By Theorem 2, in the right hand side of (\ref{NCN}),
 the terms corresponding to partitions
$\pi$ having a class containing two elements $j,k$ such that $a_j$ and $a_k$
belong to distinct algebras give a zero contribution. Thus we have to sum over
partitions in which all $j's$ belonging to a certain block of the partition
 are such that $a_j$
belongs to the same algebra. Since we can
 express noncrossing cumulants in terms of
moments we get the formula for $\tau(a_1\ldots a_n)$
 in terms of the restrictions of $\tau$ to each of the
subalgebras $A_i$. Noncrossing
 cumulants are a powerful tool for  making
computations in free probability, see \cite{NS1}, \cite{NS2}, \cite{NS3},
\cite{NSS}, \cite{SS},  for some applications.
 We give a simple   illustration  below.

Let   $X_1$ and $X_2$ be two self-adjoint elements  which are free, then the
distribution of $X_1+X_2$,  depends only on the distributions of $X_1$
 and $X_2$ and
can be computed as follows. Let $R^{(n)}(X_1,\ldots,X_1)$
and $R^{(n)}(X_2,\ldots,X_2)$, for $n\geq 1$,
 be the noncrossing cumulants of $X_1$ and $X_2$,
 then one can
expand
$R^{(n)}(X_1+X_2,\ldots, X_1+X_2)$ by multilinearity as
$\sum_{i_1,\ldots,i_n}R^{(n)}(X_{i_1},\ldots,X_{i_n})$ where the sum is over all
sequences of 1 and 2. By Theorem 2, all terms vanish except
$R^{(n)}(X_1,\ldots, X_1)$ and $R^{(n)}(X_2,\ldots, X_2)$. It follows that
$$R^{(n)}(X_1+X_2,\ldots, X_1+X_2)=R^{(n)}(X_1,\ldots, X_1)+
R^{(n)}(X_2,\ldots, X_2)$$ allowing the computation of the moments of
$X_1+X_2$, hence its distribution, in terms of the distributions of $X_1$ and
$X_2$.
It remains to give a compact form to the  relation between
moments and  noncrossing cumulants.
For any self-adjoint element $X$ with  distribution $\mu$, let
 $$G_X(z)=\frac{1}{ z}+\sum_{k=1}^\infty z^{-k-1}\tau(X^k)=\int_{\bf R}
 \frac{1}{
 z-x}\mu(dx)$$ be its Cauchy transform, and let $$K(z)=
 \frac{1}{z}+\sum_{k=0}^\infty R_kz^k$$
  be the inverse series for composition.

 {\bf Theorem 3. \cite{S1}} {\it
     $$\text{One has} \qquad R_k=R^{(k)}(X,\ldots,X)\qquad \text{ for all
     $k$.}$$}
\indent The operation which associates to the two distributions of
$X_1$ and $X_2$ the distribution of their sum is called the free
convolution of measures on the real line, and was introduced by
D. Voiculescu, who first considered the coefficients $R_k$ and
proved
 the formula  for the free
convolution of two measures,
 using
very different methods \cite{V2}.

Combining theorems 1 and 2,  given two large  random matrices of known
spectra one can predict the spectral distribution of their sum, with a good
accuracy and probability close to 1.
It is illuminating to look at the following example. The histogram below
is made
of the $800$ eigenvalues of a random matrix of the form $\Pi_1+\Pi_2$ where
$\Pi_1$ and $\Pi_2$ are two orthogonal projections onto some random subpaces of
dimension 400 in ${\bf C}^{800}$, chosen independently.
 The curve  $y=\frac{40}{ \pi\sqrt{x(2-x)}}$
which corresponds to the large $N$ limit predicted by free probability has
been drawn.
$$
\includegraphics[scale=0.5]{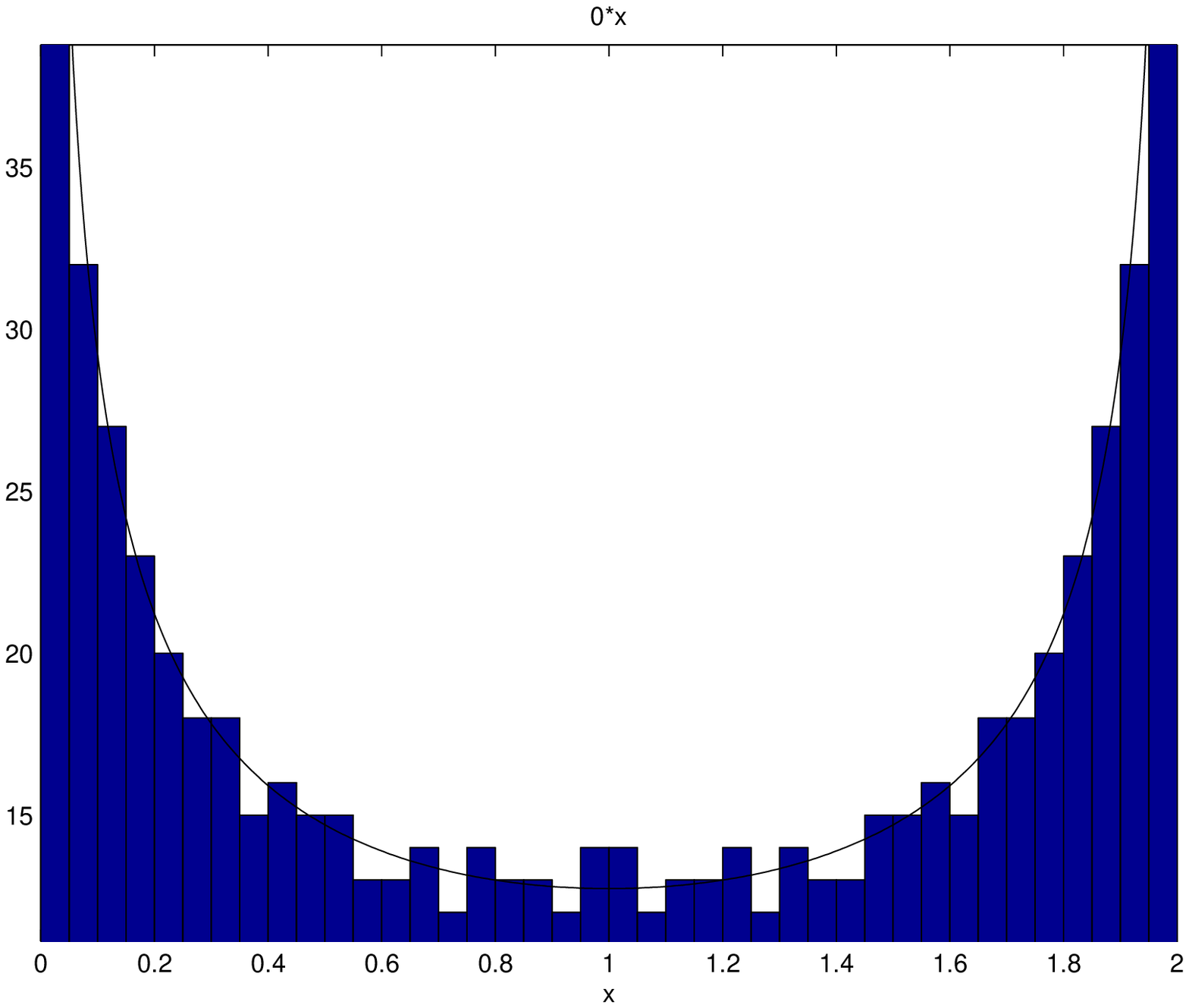}
$$

\section{Noncrossing cumulants, random matrices and characters of symmetric
groups } \label{section 4}

\setzero\vskip-5mm \hspace{5mm}

Besides free probability theory,
 noncrossing partitions appear in several areas of
 mathematics. We  indicate some relevant connections. The first is with the
 theory of map enumeration initiated by investigations of theoretical physicists
 in two-dimensional quantum field theory. The noncrossing partitions appear
 there under the guise of planar diagrams, the Feynman diagrams which dominate
 the matrix integrals in the large $N$ limit. This is of course related to the
 fact that large matrices model free probability.
  We shall not discuss this further  here,
 but refer to \cite{Z} for an accessible introduction. Another place where
 noncrossing partitions play a role, which is closely related to the preceding,
 is the geometry of the symmetric group, more precisely of its Cayley graph.
 Consider the (unoriented) graph whose vertex set is the symmetric group
 $\Sigma_n$, and such that
  $\{\sigma_1,\sigma_2\}$ is an edge if and only if $\sigma_1^{-1}\sigma_2$ is a
  transposition,
  i.e. this is the Cayley graph of $\Sigma_n$
  with respect to the generating set
  of all transpositions. The distance  on the graph is given by
  $$d(\sigma_1,\sigma_2)=n-\text{number of orbits of }\sigma_1^{-1}\sigma_2:=|
  \sigma_1^{-1}\sigma_2|.$$
  The lattice of noncrossing partitions can be imbedded in $\Sigma_n$ in the
  following way \cite{K}, given a noncrossing partition of
   $\{1,\ldots, n\}$, its image is the permutation $\sigma$
    such that $\sigma(i)$ is the
   element in the same class as $i$, which follows $i$ in the cyclic order
   $12...n$. One can check \cite{B1} that the image of $NC(n)$ is the set of all
   permutations satisfying
   $|\sigma|+|\sigma^{-1}c|=|c|$ where $c$ is the cyclic permutation $c(i)=i+1
   \, mod(n)$, in other words, this set consists of all permutations which lie on a
   geodesic from the identity to $c$ in the Cayley graph. These facts are at the
   heart of the connections between free probability, random matrices and
   symmetric groups. As an illustration we shall see how free cumulants arise
   from asymptotics of
   both random matrix theory and symmetric group representation theory.

 Recall that  cumulants (also called semi-invariants, see e.g. \cite{S})
 of a random variable $X$
 with moments of all orders,
  are
  the coefficients in the Taylor expansion of the logarithm of its
  characteristic function, i.e.
  $$\log E[e^{itX}]=\sum_{n=0}^{\infty} (it)^n\frac{C_n(X)}{ n!}.$$
  We shall consider random variables of the following form
  $Y^{(N)}=NX^{(N)}_{1,1}$ where $X^{(N)}=UD^{(N)}U^*$
  is a random matrix chosen as in (\ref{rm})
  and $X_{1,1}^{(N)}$ is its upper left coefficient.
  Assume now that the moments of $X^{(N)}$
   converge $$\frac{1}{N}Tr((X^{(N)})^k)\to_{N\to\infty}
   \int_{\bf R}x^k\mu(dx)$$ for some
  probability measure $\mu$ on $\bf R$, with
  noncrossing cumulants $R_n(\mu)$,  then one has
  $$\lim_{{N}\to\infty}\frac{1}{N^2}C_n(Y^{(N)})=\frac{1}{n}R_n(\mu).$$
 This was first observed by P. Zinn-Justin \cite{ZJ}, a proof using
 representation theory
 has been found by B. Collins \cite{C}.

 We have related noncrossing cumulants to usual cumulants
 via random matrix theory, we shall see that
  that noncrossing cumulants are also useful in evaluating characters
 of symmetric groups. The precise relation however is not obvious at first
 sight.

 Let us  recall a few facts about irreducible representations
 of  symmetric groups. It is well known that they can be
   parametrized by Young diagrams. In the following
   it will be convenient to represent a Young diagram by  a function
$\omega:\bf R\to \bf R$
 such that $\omega(x)=|x|$ for $|x|$ large enough, and
$\omega$ is a piecewise affine function, with slopes $\pm1$, see the following
picture which shows the Young diagram corresponding to the partition
$8=3+2+2+1$.
$$\includegraphics{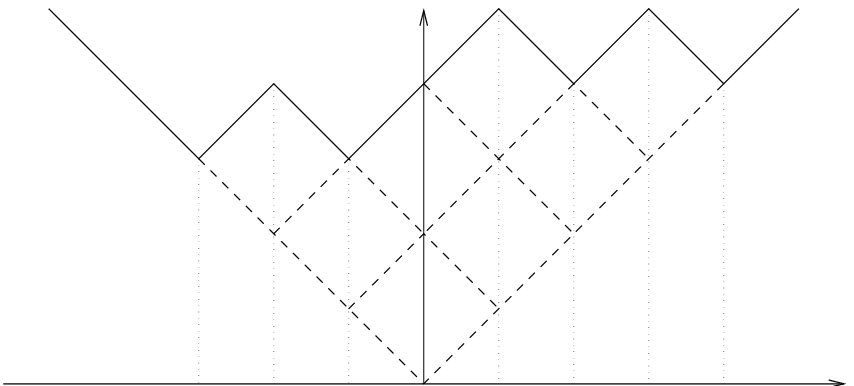}$$
  $$\qquad\qquad\ x_1\quad\,  y_1\quad
x_2\qquad\qquad  y_2\quad\, x_3\quad y_3\quad x_4\qquad$$

Alternatively we can encode the Young diagram using the local
minima and local maxima of the function $\omega$, denoted by
$x_1,\ldots, x_k$ and $y_1,\ldots, y_{k-1}$ respectively, which
form two interlacing sequences of integers. These are (-3,-1,2,4)
and (-2,1,3) respectively in the above picture.
 Associated with the Young diagram there is a unique probability
measure $m_{\omega}$ on the real line, such that
$$
\int_{\bf R}\frac{1}{
z-x}m_{\omega}(dx)=\frac{\prod_{i=1}^{k-1}(z-y_i)}{
\prod_{i=1}^{k}(z-x_i)}\quad \text{for all }z\in\bf C\setminus\bf
R.
$$
This probability measure is supported by the set $\{x_1,\ldots, x_k\}$ and is
called the transition measure of the diagram, see \cite{K1}.
 Let $\sigma$ denote the conjugacy class in $\Sigma_n$
of a permutation with $k_2 $ cycles of
length 2, $k_3$ of length 3, etc.. Here $k_2,k_3,...$ are fixed while we
let $n\to\infty$.  Denote by
 $\chi_\omega$ the normalized character of $\Sigma_n$ associated with the
Young diagram $\omega$,
 then the following asymptotic evaluation holds uniformly on the set of
 $A$-balanced Young
 diagrams, i.e. those
  whose longest row and longest column are less that $A\sqrt n$ (where $A$ is
  some constant $>0$),
\begin{equation}\label{as}\chi_\omega(\sigma)=\prod_{j=2}^{\infty}
n^{-jk_j}R_{j+1}^{k_j}(\omega)+ O(n^{-1-|\sigma|/2}).
\end{equation}
Note that $R_k$  is scaled by $\lambda^k$ if we scale the diagram $\omega$ by a
factor $\lambda$, therefore the first term in the right hand side is of order
$O(n^{\sum _j(j+1)k_j/2-\sum_j jk_j})=O(n^{-|\sigma|/2})$, this gives the order
of magnitude of the character of a fixed conjugacy group for an $A$-balanced
diagram.

 In \cite{B2}  a proof of (\ref{as}) has been given,
 using in an essential way the Jucys-Murphy
 operators. Another proof, leading to an exact formula for characters of cycles
 due to S. Kerov \cite{K2},
 was shown to me later by A. Okounkov \cite{O}, see \cite{B5}.

\section{Representations of large symmetric groups} \label{section 5}

\setzero\vskip-5mm \hspace{5mm}

 The asymptotic formula (\ref{as}) shows in particular that irreducible
 characters of
 symmetric groups become asymptotically multiplicative i.e.
 for permutations with disjoint supports $\sigma_1$ and $\sigma_2$,  one
 has \begin{equation}\label{fact}\chi_{\omega}(\sigma_1\sigma_2)=
  \chi_{\omega}(\sigma_1)\chi_{\omega}(\sigma_2)+O(
  n^{-1-|\sigma_1\sigma_2|/2})\end{equation}  uniformly on
 $A$-balanced  diagrams. Conversely, given a central, normalized,
 positive definite  function on $\Sigma_n$,  a factorization property such
 as (\ref{fact}) implies that the positive function is essentially an irreducible
  character \cite{B3}. More precisely, recall that a
  central normalized positive definite function $\psi$ on $\Sigma_n$
    is a convex combination of normalized
 characters, and as such it defines a probability measure on the set of Young
 diagrams.   For any $\varepsilon, \delta>0$, for all $n$ large enough,
 if an approximate factorization such
 as (\ref{fact}) holds for $\psi$, then there exists a curve $\omega$, such that
 the measure on Young diagrams
  associated with $\psi$ puts a mass larger than
 $1-\delta$ on Young diagrams which lie in a  neighbourhood of this curve,
  of width $\varepsilon\sqrt n$. Therefore one can say that condition
  (\ref{fact}) on a positive definite function implies that the representation
  associated with this function is approximately isotypical, i.e. almost all
  Young diagrams occuring in the decomposition have a shape close to a certain
  definite curve.

  Using this fact it is  possible to understand the asymptotic behaviour of
  several operations in representation theory. Consider for example the
  operation of induction. One starts with two irreducible
   representations of symmetric groups
  $\Sigma_{n_1},\Sigma_{n_2}$, corresponding to two Young diagrams $\omega_1$
  and $\omega_2$. One can then induce the product representation
  $\omega_1\otimes\omega_2$
  of
  $\Sigma_{n_1}\times \Sigma_{n_2}$ to $\Sigma_{n_1+n_2}$. This new
  representation is reducible and the multiplicities of irreducible
  representations
   can be computed using a combinatorial
  device, the Littlewood-Richardson rule. This rule however gives little light
  on the asymptotic behaviour of the multiplicities. Using the
  factorization-concentration
  result, one can  prove that when $n_1$ and $n_2$ are very large, but of the
  same order of magnitude, then there exists a curve, which depends on
  $\omega_1$ and $\omega_2$, and such that the typical Young diagram occuring in
  the decomposition of the induced representation, is close to this curve.
 As we saw in section 4, one can associate a probability measure on the real
 line to any Young
 diagram. The description of the typical shape of
  Young diagram which occurs in the
 decomposition of the induced representation is easier if we use this
 correspondance between probability measures and Young diagrams, indeed the
 probability measure associated with the shape of the
 typical Young diagram corresponds to the free convolution of the two
 probability measures \cite{B2}.

 There are analogous results for the restriction of representations from large
 symmetric groups to smaller ones. There the corresponding operation on
 probability measure is called the free compression, it corresponds at the level
 of the large matrix approximation, to taking a random matrix with prescribed
 eigenvalue distribution, as in section 2, and extracting a square submatrix.
 Finally there are also results for Kronecker tensor products of
 representations. Here a central role is played by the well known Kerov-Vershik
 limit shape, whose associated probability measure is the semi-circle
 distribution with density
 $\frac{1}{2\pi}\sqrt{4-x^2}$ on the interval $[-2,2]$, see
 \cite{B2}.

\label{lastpage}

\end{document}